\magnification \magstep1
\input diagrams.tex
\overfullrule=0pt
\font\tengoth=eufm10  \font\fivegoth=eufm5
\font\sevengoth=eufm7
\newfam\gothfam  \scriptscriptfont\gothfam=\fivegoth 
\textfont\gothfam=\tengoth \scriptfont\gothfam=\sevengoth

%
\font\tenbi=cmmib10  \font\fivebi=cmmib5
\font\sevenbi=cmmib7
\newfam\bifam  \scriptscriptfont\bifam=\fivebi 
\textfont\bifam=\tenbi \scriptfont\bifam=\sevenbi
\def\bi{\fam\bifam\tenbi}
\font\hd=cmbx10 scaled\magstep1
\def \Box {\hfill\hbox{}\nobreak \vrule width 1.6mm height 1.6mm
depth 0mm  \par \goodbreak \smallskip}
\def \coker {\mathop{\rm coker}}
\def \ker {\mathop{\rm ker}}
\def \deg  {\mathop{\rm deg}}

\def \Im  {\mathop{\rm Im}}
\def \rank {\mathop{\rm rank}}
\def \iso {\cong}
\def \tensor {\otimes}

\def \Tor {{\rm Tor}}
\def \Ext {{\rm Ext}}

\def \th {{^{\rm th}}}

\def \P {{\bf P}}
\def \H {{\rm H}}

\def \O {{\cal O}}
\def \I {{\cal I}}
\def \L {{\cal L}}
\def \N {{\cal N}}
\def \E  {{\cal E}}
\def \Sym {{\rm Sym}}
\def \h {{\rm h}}

\newarrow{Equals}=====

\def\boxit#1#2{\vbox{\hrule\hbox{\vrule
\vbox spread#1{\vfil\hbox spread#1{\hfil#2\hfil}\vfil}%
\vrule}\hrule}} 

\rightline{December 30, 1996}
\bigskip

\centerline{\hd Gale Duality and Free Resolutions}
\smallskip
\centerline{\hd of Ideals of Points}
\medskip 
\centerline {by} 
\smallskip 
\centerline {\bf David Eisenbud and Sorin Popescu
\footnote{$^{*}$}{\rm The first author
is grateful to the NSF and the second author to the DFG for 
support during the preparation of this work.}}

\bigskip
\bigskip

 What is the shape of the free resolution of the ideal of a general
set of points in $\P^r$? This question is central 
to the programme of connecting the geometry of point sets in
projective space with the structure of the free resolutions of 
their ideals.  There is a 
lower bound for the resolution computable from the (known)
Hilbert function, and it seemed natural to conjecture that
this lower bound would be achieved.  This is the
``Minimal Resolution Conjecture'' (Lorenzini [1987], [1993]).
Although the conjecture has been shown to hold in many cases,
three examples discovered
computationally by Frank-Olaf Schreyer in 1993
show that it fails in general.
In this paper we shall describe a novel structure
inside the free resolution of a set of points which accounts
for the failure and provides a counterexample in $\P^r$ for
every $r\geq 6, \ r\neq 9$.

We begin by reviewing the conjecture and its status.
Consider a set of $\gamma$ points in the projective $r$-space
over a field $k$, say $\Gamma\subset\P_k^r$.
Let $S=k[x_0,\dots,x_r]$, let $I_{\Gamma}$ be the homogeneous ideal of 
$\Gamma$, let $S_{\Gamma}$ denote the homogeneous coordinate ring 
of $\Gamma$. Let 
$$
F_\bullet:\quad 0 \rTo F_{r-1} \rTo \ldots 
\rTo  F_0 \rTo I_{\Gamma} \rTo 0
$$
be the minimal free resolution of $I_{\Gamma}$, and define the 
associated (graded) betti numbers $\beta_{ij}$ by the formula
$$
F_i = \oplus_j S(-j)^{\beta_{ij}}.
$$
The minimal free resolution conjecture can be formulated as follows:

\proclaim Minimal Resolution Conjecture. If $\Gamma$ is a general
set of points in $\P_k^r$ over an infinite field $k$, then for
any integers $i,j$, at most one of
$\beta_{i,j}$ and $\beta_{i+1,j}$ is nonzero.

Given our knowledge of the Hilbert function of 
the general set of points (since $\Gamma$ imposes independent 
conditions on forms of every degree) and the easy result that if 
$I_{\Gamma}$ contains forms of degree $d$, then $\beta_{ij}=0$
for $j>i+d$ -- that is, $I_{\Gamma}$ is $(d+1)$-regular,
the minimal free resolution conjecture can 
be translated into an explicit formula for the $\beta_{ij}$
(see \S 5 below).

The minimal resolution conjecture is known to be true in $\P^2$
(Gaeta [1951] and [1995], Geramita-Lorenzini [1989]), in 
$\P^3$ (Ballico-Geramita [1986]), in $\P^4$
(Walter [1995], Lauze [1996]), and in $\P^n$ for
$n+1\le\gamma\le n+4$, or 
$\gamma={{n+2}\choose 2}-n$ (Geramita-Lorenzini [1989],
Cavaliere-Rossi-Valla [1991], Lorenzini [1993]).
Its predictions about $\beta_{r-1,j}$
are known to be true in general (Lauze [1995]).  
Most striking, the conjecture is known to hold whenever 
the number of points in $\Gamma$ is sufficiently large compared to 
$r$ (Hirschowitz-Simpson [1994]), where the bound 
given is $\gamma>6^{r^3\log r}$.  

Schreyer discovered by computational experiments of a probabilistic
nature that the following three cases give counterexamples to the
conjecture: 11 points in $\P^6$, 12 points
in $\P^7$, and 13 points in $\P^8$. A few more such examples
were discovered by computer search (Boij [1994],
Beck-Kreuzer [1996]). Despite considerable
effort, no-one was able to give a 
non computational treatment of these
examples, nor to find any ``explanation'' of them,
so that it was unclear whether they were unique accidents
or part of a larger picture. 

In this paper we give a geometric construction that gives
rise to a subcomplex of the resolution of a general set of points.
A consequence of our construction is the following, which includes
all the examples that are currently known:

\proclaim Theorem 0.1. For any integer $r\geq 6,\  r\neq 9$, there is
an integer $\gamma(r)$ such that the Minimal Resolution Conjecture
fails for a set of $\gamma(r)$ general points in $\P^r$. 
More explicitly, if $s$ and $k$ are (uniquely) defined by
$$r={s+1\choose 2}+k,\quad 0\le k \le s,$$
then we may take 
$$\gamma(r) = r+s+2 = {s+2\choose 2}+k+1.$$

We do not know whether such examples exist in $\P^9$; but Beck-Kreuzer
[1996] have made computations showing that
none occurs for 50 or fewer points.
\medskip

Here is an outline of the ideas involved.  Associated to
any embedding $\Gamma\subset\P^r$ of a set of $\gamma$ points 
(in sufficiently general position) 
in projective space is another embedding of the
same set of points in projective $(s:= \gamma -r-2)$-space, called the
{\it Gale Transform} of $\Gamma$ (see \S 1).
Call the image of the transformed 
embedding $\Gamma'\subset\P^s$; it is again a general set of points.  
We may identify the ambient space $\P^r$ of $\Gamma$ with the space of
lines in $\H^1(\I_{\Gamma'}(1))$.
Using this identification, we relate the back ends of the resolutions
of $\Gamma$ and $\Gamma'$ (see \S 1). 
Writing $W := \H^0(\O_{\P^s}(1))$, and $U :=\H^1(\I_{\Gamma'}(2))^\ast$
we have a multiplication pairing
$$
\mu:\quad W\tensor U\rTo \H^1(\I_{\Gamma'}(1))^\ast=\H^0(\O_{\P^r}(1)).
$$
Associated to any such pairing is
a complex built from a certain Koszul complex 
$$
E_\bullet(\mu):\ \ldots \rTo\wedge^3 W\tensor D_2 U\tensor\O_{\P^r}(-2)  \rTo 
\wedge^2 W\tensor U\tensor\O_{\P^r}(-1)  \rTo W\tensor\O_{\P^r},
$$
where $D_lU$ denotes 
the $l^{\rm th}$ divided power of $U$ (see \S 3).
There is a natural map  ${E_\bullet(\mu)(r+2)}$ 
into the dual $F_\bullet^\ast$ of the free resolution of $\I_\Gamma$, 
regarded as a complex of sheaves (we sometimes regard $E_\bullet(\mu)$ as
a complex of free modules).
When the map $\mu$ satisfies a certain nondegeneracy 
condition, the complex $E_\bullet(\mu)$ has a property we call
linear exactness (see \S 2). In this situation
the map  $E_\bullet(\mu)(r+2)\rTo F_\bullet^\ast$ is a monomorphism
onto a direct summand, and this gives a lower bound for the 
betti numbers $\beta_{ij}$ of $\I_\Gamma$ that is sometimes in conflict
with the conclusion of the Minimal Resolution Conjecture. One
might say, in summary, that the failure of the conjecture for
a set $\Gamma$ comes from the failure
of $\Gamma'$ to impose independent conditions on forms
of degree 2; 
not because $\Gamma'$ isn't sufficiently 
general, but because its degree is greater than the number of
forms of degree 2 in $\P^s$.

The  heart of the paper, and by far
the most difficult part,  is 
the nondegeneracy of the pairing $\mu$, established in \S 4.  
This nondegeneracy represents an open condition
on the family of sets of points $\Gamma$.  Thus in order to prove
that it holds for the general $\Gamma$, it is enough to show that
it holds for a  special set of points $\Gamma$. We do this
by specializing the points to lie on a curve $C$. Under favorable
circumstances, the nondegeneracy condition on $\mu$ can 
be re-interpreted as a cohomology condition on a certain vector
bundle on the curve $C$; the argument has the flavor of Koszul
cohomology.  We complete the argument by specializing $C$
to either a plane curve or to a  curve with prescribed gonality
(depending on the parity of $s$), and taking the points in such 
a way that the bundle in question decomposes into a direct sum of 
simpler bundles; even then the computation of cohomology involves
some nonstandard ideas.  For instance, in the plane curve case
case $s=4$, we must show the following: Let $C$ be
a general plane curve of degree $\geq 5$ and let $T$ be the restriction of
the tangent bundle of the plane to $C$. If $L$ is a  general
line bundle on $C$ of degree genus$(C)-1$, then there are
no ``twisted endomorphisms'' $T \rTo L\otimes T$.
We can actually prove the vanishing
theorem in the plane curve case (Theorem 4.1) only in characteristic 0;
but as Theorem 0.1 in full generality follows from the
case of characteristic $0$, this is no problem.

To see how all this works in the easiest interesting case,
let $\Gamma$ be a set of 
$\gamma(6)=11$ general points in $\P^6$, Schreyer's simplest example.
With notation as in Theorem 0.1 we have $s$ = 3.
If we display the betti numbers $\beta_{ij}$ associated
with the resolution $F_\bullet$ in a table in the style 
of the program Macaulay,
the expected betti numbers, coming from the Minimal 
Resolution Conjecture, would be
$$
\vbox{\offinterlineskip 
\halign{\strut\hfil# \ \vrule\quad&# \ &# \ &# \ &# \ &# \ &# \ 
&# \ &# \ &# \ &# \ &# \ 
\cr
degree&\cr
\noalign {\hrule}
0&1&--&--&--&--&--&--\cr
1&--&17&46&45&4&--&--\cr
2&--&--&--&--&25&18&4\cr
\noalign{\bigskip}
\omit&\multispan{8}{\bf Conjectural shape of $F_\bullet$}\cr
\noalign{\smallskip}
}}
$$
We have $\dim(W) = 4$.
The set $\Gamma'$ consists of 11 general points in $\P^3 = \P(W)$,
and since the space of quadrics in $\P^3$ is only 10-dimensional,
$\dim(U) = 1$. Identifying the
divided powers of $U$ with the ground field, the complex 
$
E_\bullet(\mu)
$
becomes
$$
E_\bullet(\mu):\ 0 \rTo\wedge^4 W\tensor\O_{\P^r}(-3) 
\rTo\wedge^3 W\tensor\O_{\P^r}(-2)  \rTo 
\wedge^2 W\tensor\O_{\P^r}(-1)  \rTo W\tensor\O_{\P^r},
$$
which can be identified with the back end of (a twist of)
the Koszul complex of a sequence of four linear forms on $\P^6$.
The nondegeneracy condition on $\mu$ becomes
the condition that the complex $E_\bullet(\mu)$ is exact (in this
case exactness and linear exactness coincide). 
The nondegeneracy can be proved by a special
argument in this case (see below),
but our general method is the following:  Since the condition
on $\Gamma$ (or, equivalently, on $\Gamma'$) is open, it suffices to
prove the result after degenerating 
$\Gamma'$ until it lies on a curve $C\subset \P^s = \P^3$.  In this case
we take $C$ to be a general sextic curve in $\P^3$ of genus 3,
and let $H$ denote its hyperplane class. Such a curve is 
projectively normal.  By a 
Koszul homology argument we show that the nondegeneracy condition
we need follows if we can show that 
$$
\H^0(\wedge^2 M_H\tensor\O_C(K_C+\Gamma^\prime-2H))=0,
$$
where $M_H$ denotes the rank 3 vector bundle that is
the kernel of the evaluation map 
$\H^0(\O_C(H))\tensor \O_C \rTo \O_C(H)$, and $K_C$
is the canonical class of $C$.  We may think of the line bundle
$\O_C(L) :=\O_C(K_C+\Gamma^\prime-2H)$ simply as a general
line bundle of degree $(2g(C)-2) +11- 2\cdot 6 = 3$ on $C$. 
To  prove the required vanishing, we degenerate $C$ to a curve
of type $(2,4)$ on a smooth quadric surface in $\P^3$
(so that $C$ is hyperelliptic), and  $\O_C(H)$ to $\O_C(3H_0)$,
where $\O_C(H_0)$, the line bundle corresponding to the hyperelliptic
involution on $C$, is induced by the class  $(0,1)$
on the quadric.  The canonical series on $C$ is induced by
$(0,2)$ on the quadric, from which we easily see that 
$\O_C(3H_0)$
is nonspecial and that
$$
M_{3H_0}\iso \oplus_{i=1}^3 \O_C(-H_0)
$$ 
is a degeneration
of $M_H$.  Thus it suffices to show that 
$$
\H^0((\wedge^2 M_{3H_0})\tensor \O_C(L)) =
\oplus_{i=1}^3 \H^0(\O_C(L-2H_0)) =0.
$$
But $\O_C(L-2H_0)$ is a  line bundle of degree
$3-4= -1$, so the result is now immediate. 
The same argument works for 12 or 13 points on $C$, giving the cases in 
$\P^7$ and $\P^8$ respectively. 
In the cases $r=6$ and $r=7$ (but already not for $r=8$)
the necessary nondegeneracy can be proved more simply: the pairing
$\mu$ can be identified as the multiplication map of sections of 
certain line bundles on the curve $C$ (see \S 4) and as such
is $1$-generic (in the sense of Eisenbud [1988]). 
In general, the method works
when $\dim U\le 2$, and  Kreuzer [1994] has proved
the necessary $1$-genericity in all cases, but the nondegeneracy
we need does not follow from this when $\dim U> 2$.

As shown in the body of the paper, it follows that the 
complex 
$
E_\bullet(\mu)(8)
$
is a subcomplex of the dual $F_\bullet^\ast$ of the minimal free resolution
of $\I_\Gamma$. Equivalently, $F$ maps onto the complex 
$
E_\bullet(\mu)(-8)^\ast
$
(suitably shifted) which has betti
display
$$
\vbox{\offinterlineskip 
\halign{\strut\hfil# \ \vrule\quad &# \ &# \ &# \ &# \ &# \ &# \ &# \ \cr
degree&\cr
\noalign {\hrule}
0&--&--&--&--&--&--&--\cr
1&--&--&--&--&--&--&--\cr
2&--&--&--&1&4&6&4\cr
\noalign{\bigskip}
\omit&\multispan{5}{$E_\bullet(\mu)(-8)^*$}\cr
}}
$$

Under these circumstances {\it each\/} betti number for $F_\bullet$
must be at least as big as the one for $E_\bullet(\mu)(-8)^\ast$, 
so we see that a lower bound for the size of $F_\bullet$ is given
by the following betti diagram (we have indicated
the differences inside boxes). In this case $r=6$ (and also when
$r=7$, but already not when $r=8$), computation shows 
that this diagram gives the actual value of the $\beta_{ij}$, so that
the theory here developed leads to an exact computation.
$$
\vbox{\offinterlineskip 
\halign{\strut\hfil# \ \vrule\quad&# \ &# \ &# \ &# \ &# \ &# \ 
&# \ &# \ &# \ &# \ &# \ &# \ &# \
\cr
degree&\cr
\noalign {\hrule}
0&1&--&--&--&--&--&--&\cr
1&--&17&46&45&\boxit{3pt}{5}&--&--&\cr
2&--&--&--&\boxit{3pt}{1}&25&18&4&\cr
\noalign{\bigskip}
\omit&\multispan{6}{\bf Actual shape of $F_\bullet$}\cr
\noalign{\smallskip}
}}
$$

The four linear forms that enter the
pairing $\mu$ in the case $r=6$ have an amusing
interpretation: They define a plane $\Pi$ in $\P^6$,
which is distinguished just by the data of the 11 points;
it is defined from the structure of the
cohomology module of the ideal sheaf of the points in $\P^3$.
The plane $\Pi\subset\P^6$ is spanned by (any) three points which
together with the 11 initial ones form a collection of 
self-associated points in $\P^6$ (that is a set which is
self-dual with respect to the Gale Transform).
Charles Walter has pointed out to
us that $\Pi$ could also be interpreted as the unique plane in $\P^6$ such
that the projection of the 11 points from this plane into $\P^3$ is
equivalent to the the Gale transform of the 11 points. 
(This latter characterization follows directly from the theory in \S 1.)

It is interesting to compare the case of points with that of curves.
The  minimal resolution conjecture for complete embeddings
of large degree (compared with the genus) general curves was shown
to be false by Schreyer [1983], and Green-Lazarsfeld [1988];
the failure comes essentially from the existence of special
divisors on the curve, which give rise to rational normal scrolls
containing the curve, and is quite different in character from the
phenomena exhibited here.  
By contrast, no counterexamples to the appropriate
minimal free resolution conjecture are known for ideals in a polynomial
ring which are made from a generic vector space of forms of some
degree $d$ plus all the forms of degree $d+1$; these are the 
ideals that seem to be the most reasonable analogue of
ideals of general sets of points.  However, the problem is
computationally difficult, and not many cases have been 
examined.

It is a pleasure to thank Mike Stillman, who joined us in
discussions leading
to some of the ideas in this paper, Andr\'e Hirschowitz and
Charles Walter, from whose ideas the exposition has benefitted, and  
Bob Friedman, who pointed out to us the beautiful paper of Raynaud [1982].
We are also
grateful to Stillman and to Dave Bayer for the program {\it Macaulay\/}
(Bayer-Stillman [1989--1996]) which has been extremely useful to us; without
it we would probably have never been bold enough to guess the
existence of the structure that we explain here. Finally, 
we are grateful to Mark Green: in earlier
joint work the first author learned from him how useful
maps on the cohomology of the ideal sheaves of points could be;
this helped to spot the connection exploited in this paper.

\beginsection \S 1. The Gale Transform

We first give a naive definition
of the Gale transform of a set of points.
Then we explain a more flexible view, in which the Gale transform
is an involution --- essentially Serre duality --- on the set of 
linear series on a set of points.  
Finally, we exhibit a peculiar module which maps to the
canonical module of a suitable set of points.  This module 
has a natural interpretation in terms of the Gale transform.  In the
next section we will exhibit a subcomplex of the resolution of this module
that is ``responsible'' for the failure of the minimal resolution conjecture.
\medskip

\noindent{\sl Definition.\ }  Let $k$ be a field, and write
$\P^r$ for $\P^r_k$.
Let $\Gamma \subset \P^r$ be a set of $\gamma$ labelled points such
that every subset of $\gamma -1$ of the points
spans $\P^r$. Choosing
homogeneous coordinates for the points, we may write their
coordinates in the form of a matrix $G:\ k^{r+1}\rTo k^\Gamma$, and
this matrix has rank $r+1$.  If
we dualize this matrix and take the kernel, we get a matrix
$G':\ k^{s+1}\rTo (k^\Gamma)^*$, where $s+r+2 = \gamma$. 
Since $k^\Gamma$
has (up to scalars) a natural basis, consisting of functions
vanishing at all but one point, we may identify 
$k^\Gamma$ and $(k^\Gamma)^*$ in a way that is natural up to the
choice of a diagonal matrix, and regard $G'$ as being a map
$k^{s+1}\rTo k^\Gamma$.  The rows of this matrix determine points
in a set $\Gamma'$ of labelled points
in $\P^s$, labelled by the same set as $\Gamma$; the rows are
all nonzero because of the condition that 
every subset of $\gamma -1$ of the points
spans $\P^r$.
The reader may check that $\Gamma'$ is
uniquely determined from $\Gamma \subset \P^r$ up to the action
of $PGL(s+1)$, and that it spans $\P^s$.  
The set $ \Gamma'$ is called the (classical) {\it Gale transform\/}
of $\Gamma$. 

\medskip
 
The Gale transform has a long history, extending at least
as far as Hesse's thesis [1840]. It was studied 
(under the rubric ``associated sets of points'')
by Castelnuovo [1889] 
and later by A.B.~Coble ([1915, 1916, 1917, 1922]), 
who discovered amazing geometric constructions
 and applied the Gale transform
to the study of Theta-functions and Jacobians
in the early part of this century. For a modern exposition
with many extensions see Dolgachev and Ortland 
[1988].  The Gale transform has reappeared in many places.
For example V.~D.~Goppa reinvented the idea in 
 the context of coding theory.  He proved
that  if $\Gamma$ lies on a linearly normal
curve $C\subset\P^r$, then $\Gamma^\prime$ lies on a different 
embedding of the same curve (see Goppa, [1984]).
The name ``Gale transform'' has become established
by  the very fruitful use
of the idea (in a somewhat different form) in the study of
convex polytopes and integer programming initiated by
D.~Gale in [1963].
We refer to the forthcoming paper Eisenbud-Popescu [1997]
for more history and geometric constructions. 
\medskip

Recall that a linear series on a scheme $X$ is a pair $(V,L)$ consisting
of a line bundle $L$ and a vector space $V$ of global sections of $L$.
The Gale transform can be defined much more generally, as
an involution on the space of 
linear series on a finite scheme $\Gamma$.  
Of course it is somewhat pedantic to speak of line bundles and
global sections on
a finite scheme, since any such scheme is affine and every line bundle
is trivial, but it has the same virtues as does the distinction
between a vector space and its dual:
this language will allow us to make definitions without
any arbitrary choices.

If $\Gamma$ is a Gorenstein scheme, finite over a field $k$, 
and $L$ is a line bundle on $\Gamma$, then
Serre duality provides a canonical ``trace''
$\tau:\ \H^0(K_{\Gamma})\rTo k$ with the property that for any line bundle
$L$ on $\Gamma$ the composition 
$$
\H^0(L)\tensor_k \H^0(K_{\Gamma}\otimes L^{-1}) \rTo \H^0(K_{\Gamma}) \rTo k.
$$
of $\tau$ with the multiplication map gives a perfect pairing  between
$\H^0(L)$ and $\H^0(K_{\Gamma}\otimes L^{-1})$.
If $V\subset \H^0(L) $ is a subspace, then we write 
$V^\perp \subset \H^0(K_{\Gamma}\otimes L^{-1})$ for the annihilator.

Using these idea, we may define the Gale transform more generally:
\smallskip
\noindent{\sl Definition.\ }  Let $k$ be a field, and let
$\Gamma$ be a Gorenstein scheme finite over $k$.
The {\it Gale transform} of a linear series $(V,L)$  on $\Gamma$  
is the linear series
$(V^\perp,K_{\Gamma}\otimes L^{-1})$. (This is the natural definition of 
adjoint series in the zero-dimensional case.)

We recall that the ``Veronese'' linear
series are defined by multiplication: 
If $n\geq 1$, then we write $V^n$ for the image by multiplication
of $V^{\tensor n}$ in $\H^0(L^n)$.  We set
$V^0 = k\subset \H^0(\O_\Gamma)$, while for $n<0$ we set $V^n = 0$ 
(again as a subset of $\H^0(L^n)$).  

The relation to the classical Gale transform is included in the
following alternative description:

\proclaim Proposition 1.1. Let $k$ be a field, and let
$\Gamma$ be a Gorenstein scheme finite
over $k$.
If $r\ge 1$ and the linear series
$(V,L)$ defines an embedding of $\Gamma$ in
$\P^r=\P(V)$ with ideal sheaf $\I_\Gamma$,
then there are natural identifications 
$(V^n)^\perp = \H^1(\I_\Gamma(n))^\ast$.
If further $\Gamma$ is a reduced set of $k$-rational
points and every subset of $\gamma -1$ of the points of 
$\Gamma$
spans $\P^r$, then the linear series
$(V^\perp,K_{\Gamma}\otimes L^{-1})$ is base-point-free,
and  the image of
$\Gamma$ under the corresponding map is the classical 
Gale transform of $\Gamma$.

\noindent{\sl Proof.\ } Using Serre duality to identify 
$\H^0(K_{\Gamma}\otimes L^{-n})$ with the dual of $\H^0(L^n)$,
the space $(V^n)^\perp$ becomes the kernel of the map
$\H^0(L^n)^\ast \rTo (V^n)^*$ dual to the inclusion.  In
the setting of the classical Gale transform 
we choose an identification of
$\H^0(L)$ with $k^\Gamma$, and the last statement of the Proposition
follows.  More generally, if $(V,L)$ defines an embedding 
of $\Gamma$,
then the exact sequence 
$$0\rTo \I_\Gamma \rTo \O_{\P(V)} \rTo \O_\Gamma \rTo 0$$
gives rise to an exact sequence 
$$
\bigl[\Sym_nV= \H^0(\O_{\P(V)}(n))\bigr]
\rTo 
\H^0(L^n)
\rTo \H^1(\I_\Gamma (n))
\rTo \bigl[0 = \H^1(\O_{\P(V)}(n))\bigr],
$$
which yields the identification ${((V^n)^\perp)}^\ast = 
\H^1(\I_\Gamma(n))$ as required.
\Box

The next result gives a description of the $(V^n)^\perp$ that
does not depend on the points being embedded:

\proclaim Proposition 1.2.  Suppose that $\Gamma$ is a 
Gorenstein scheme, finite over $k$, and let
$(V, L)$ be a linear series on $\Gamma$. 
For each integer $n$ the product
$V\cdot (V^n)^\perp$ lies in $(V^{n-1})^\perp$.  If $n\neq 0$, 
then $(V^n)^\perp$ is the 
largest subspace of $\H^0(K_\Gamma\tensor L^{-n})$ that
is multiplied by $V$ into $(V^{n-1})^\perp$.

\noindent{\sl Proof.\ } If $n\leq 0$ the result is trivial.
If $n>0$ and $a \in\H^0(K_\Gamma\tensor L^{-n})$,
then $a \in (V^n)^\perp$ iff $\tau(a\cdot V^n)=0$ 
iff $\tau(aV\cdot V^{n-1})=0$
iff $aV\subset (V^{n-1})^\perp$.
\Box
\medskip

By virtue of Proposition 1.2 we may regard 
$\oplus_{n\in {\bf Z}} ({(V^{-n})}^\perp)$
as a graded $k[V]$-module, and with this structure we will call it
$\omega_{\Gamma,V}$.  In the case
where $\Gamma$ is embedded in $\P(V)$, the following Corollary
of Proposition 1.1 identifies this module with the canonical module of the 
affine cone over $\Gamma$.  Since the minimal free resolution of 
this canonical module
is the dual of the minimal free resolution of $k[V]/I_\Gamma$, this result
provides the link with free resolutions:

\proclaim Corollary 1.3.  Let
$\Gamma \subset\P(V)\iso\P^r$
be a zero-dimensional Gorenstein
subscheme, and let $S = k[V]$ be the polynomial ring
in $r+1$ variables.  Write $I_{\Gamma}$ for the homogeneous ideal
of $\Gamma$. There is a natural isomorphism
$$
\omega_{\Gamma,V} \iso \Ext^{r-1}_S(I_\Gamma, S(-r-1)).
$$

\noindent{\sl Proof.\ } By Serre duality, 
$\Ext^{r-1}_S(I_\Gamma, S(-r-1)) = \oplus_n(\H^1(\I_\Gamma(n))^*)$ as
$S$-modules, the multiplication 
$$
V \tensor \H^1(\I_\Gamma(n))^\ast
\rTo 
\H^1(\I_\Gamma(n-1))^\ast
$$
being the one induced by the multiplication
$\H^0(\O_{\P(V)}(1)) \tensor \H^1(\I_\Gamma(n-1)) \rTo \H^1(\I_\Gamma(n))$.
Proposition 1.1 identifies $\H^1(\I_\Gamma(n))^*$ with $(V^n)^\perp$ for
$n>0$, while the identification for $n\leq 0$ is trivial.  The compatibility
of these identifications with the multiplication maps follows from the
same exact sequence as employed in the proof of Proposition 1.1.
\Box

\medskip

We now approach the fundamental construction to be studied in this paper.
We write  $S_{a,b}(W)$ for the Schur functor
$$\eqalign{
S_{a,b}(W) &:= \Im (\wedge^{a+1}W\tensor \Sym_{b-1}W \rTo
\wedge^{a}W\tensor \Sym_{b}W)\cr
& =\ker(\wedge^{a}W\tensor \Sym_{b}W\rTo
\wedge^{a-1}W\tensor \Sym_{b+1}W).}
$$
For example, $S_{1,1}(W)=\wedge^2 W$, the inclusion into 
$\wedge^{a}W\tensor \Sym_{b}W = W\tensor W$ being the 
diagonal map of the exterior algebra.  
The reader unfamiliar with Schur functors may avoid them at first
by considering only this case, that is, taking $n=2$ in the following
result.

\proclaim Theorem 1.4.  Let $\Gamma$ 
be a Gorenstein scheme, finite over $k$. Let $(V,L)$ and 
$(W = V^\perp, K_\Gamma\tensor L^{-1})$ be dual linear series, and set
$U := (W^n)^\perp$, with $n\ge 1$. The natural multiplication 
 $\Sym_{n-1}W\tensor U\rTo W^{n-1}U\subset V$
induces a map 
$\mu:\ \Sym_{n-1}W\tensor U \rTo V$
which in turn defines a map of free $k[V]$-modules
$$
\delta:\quad S_{1,n-1}(W)\tensor U\tensor k[V] 
\rTo W\tensor k[V](1).
$$
There is a unique map of $k[V]$-modules
$
(\coker \delta) \rTo  \omega_{\Gamma, V}
$
which extends the inclusion 
$W = (\omega_{\Gamma, V})_{-1} \subset \omega_{\Gamma, V}.$

\noindent{\sl Proof.\/} The map $\delta$ is the composite
$$
S_{1,n-1}(W)\tensor U\tensor k[V] 
\rTo
W\tensor \Sym_{n-1}W \tensor U\tensor k[V]
\rTo^{W\tensor \mu}
W\tensor k[V](1),
$$
whereas the space of linear relations on 
$\omega_{\Gamma, V}$ can be identified
with the vector space $N$ which is the kernel of the multiplication map
$m:\ W\tensor V \rTo \H^0(K_{\Gamma})$. We must show that $N$ contains
the relations on $\coker \delta$, which are generated by the image
of the composite 
$$
S_{1,n-1}(W)\tensor U\rTo
W\tensor \Sym_{n-1}W \tensor U
\rTo^{W\tensor\mu}
W\tensor V.
$$
Now the natural multiplication map 
$\mu^\prime:\ \Sym_n W\tensor U \rTo \H^0(K_\Gamma)$ fits in the diagram
$$
\diagram[midshaft]
0&\rTo & N&\rTo&W\tensor V&\rTo^m&W\cdot V&\rInto& \H^0(K_\Gamma)\\ 
&&\uDashto&&\uTo_{W\tensor\mu}&&&&\uTo_{\mu^\prime}\\
0&\rTo&S_{1,n-1}(W)\tensor U&\rTo&W\tensor\Sym_{n-1}(W)\tensor U&&\rTo&&
\Sym_{n}(W)\tensor U,\\
\enddiagram
$$
which is commutative by the associativity of multiplication, and has exact
rows by the definition of $N$ and $S_{1,n-1}$.  Thus there is a vertical
map induced on the left, which is the desired inclusion.
\Box
\smallskip

The significance of this result is that it gives a map of 
complexes from a complex $E^{-1}_\bullet(\mu)(r+2)$ (described in \S 3)
beginning with the map $\delta$
into the dual of the resolution of the ideal of the points. We shall see
that under certain circumstances this map is an inclusion, and provides 
the subcomplex which ``spoils'' the Minimal Resolution Conjecture.
The properties of this map will be the subject of \S 2. 
Of course Theorem 1.4 is vacuous if $U = (W^n)^\perp = 0$. By Proposition 1.1,
if $\Gamma$ is a set of points in $\P^r = \P(V)$ and $\Gamma'$ is its Gale
transform, embedded in $\P^s = \P(W)$, then $U =
\H^1(\I_{\Gamma'}(n))^\ast$; 
thus $U$ is nonzero iff $\Gamma'$ fails to 
impose independent conditions on forms
of degree $n$.  The analysis of the resolution of $I_\Gamma$
via the map $\delta$ will involve the geometry of the Gale 
transform $\Gamma'$.
\medskip

\noindent{\sl Remark.\ }  There is a less invariant version of these ideas
which is pleasingly direct:
Again let $\Gamma$ be a Gorenstein scheme, finite over a field $k$,  and
let $\O_\Gamma$ be the coordinate ring of $\Gamma$, a finite dimensional
Gorenstein $k$-algebra.   Suppose that
$\Gamma$ is embedded in $\P^r$. 
If we choose a hyperplane not meeting $\Gamma$, we may identify
the line bundle $L = \O_\Gamma(1)$
with $\O_\Gamma$, and thus identify the linear series 
$(V = \H^0(\O_{\P^r}(1)), L)$ with a subspace $V\subset \O_\Gamma$. 
We also choose an identification of $\O_\Gamma$
with $K_\Gamma$ (equivalently, we may choose a ``trace'' functional
$\tau: \O_\Gamma \rTo k$ not vanishing on any component of the socle of 
$\O_\Gamma$) and consider the pairing on $\O_\Gamma$ defined as the 
composition of this functional with multiplication. We may again
define the powers $V^n$ and the spaces 
$W_{n}:= {(V^{-n})}^\perp$, but this
time they will all be subspaces of $\O_\Gamma$.

\beginsection \S 2. Linear exactness and linear rigidity

We shall use a property of certain complexes that
we have not seen exploited before: we call it
{\it linear exactness}. We give its definition in
an abstract setting before plunging into the multilinear
algebra necessary to define the complexes to which we will
apply it.

Let $S$ be a graded ring with $S_0 = k$ a field, and let 
$$
E_\bullet: \quad \ldots \rTo E_{i+1} \rTo E_i \rTo \dots \rTo E_n
$$ 
be a {\it linear complex\/} 
in the sense that each $E_i$ is a free module generated in
degree $i$, so that in particular
the differentials are given by matrices of elements of $S_1$.
  We shall say that $E$ is {\it linearly exact\/} if, for all $i>n$,
the homology $\H_i(E_\bullet)$ is nonzero only in degrees $>i$,
or equivalently if,
in any matrix representing a differential
of $E_\bullet$,  the columns are linearly independent over $S_0$. 
The utility of this
definition lies in the following result:

\proclaim Proposition 2.1. Let $E_\bullet$ be a linearly exact 
linear complex as above, and
and let 
$$
F_\bullet: \quad \ldots \rTo F_{i+1} \rTo F_i \rTo\ldots\rTo F_n
$$ 
be a graded minimal free resolution with $F_n$ generated in 
degrees $\geq n$. If
$\alpha:\ E_\bullet\rTo F_\bullet$ is a map of complexes such that
$\alpha_n:\ E_n\rTo F_n$ is an inclusion, then
each map $\alpha_i:\ E_i\rTo F_i$ is a split inclusion.

\noindent{\sl Proof.\ } Any set of 
elements of degree $n$ in $F_n$ that are linearly independent over
$S_0$ are part of a free basis. Thus
$E_n$ maps to a direct summand of $F_n$.
It follows from the hypothesis that for each $i$ 
the module
$F_i$ is generated in degrees $\geq i$. 
 By hypothesis, 
a free basis of $E_{n+1}$ maps to a set of linearly independent
elements in $(E_n)_{n+1}$, which is a summand of $(F_n)_{n+1}$. Since
$F_{n+1}$ is generated in degrees $\geq n+1$,
the map $\alpha_{n+1}$ must take the basis of  $E_{n+1}$
to a subset of a basis of $F_{n+1}$. Thus $\alpha_{n+1}$ is a split 
inclusion, and induction completes now the proof.\Box

It often suffices to prove linear exactness
at the first step:

\proclaim Lemma 2.2.(Linear rigidity). Let $R=k[x_0,\dots,x_r]$
be a polynomial ring, and let $M$ be an $R$-module generated in
degrees $\geq 0$.  Let $F_\bullet$ as above be a minimal free resolution of $M$.
Suppose that $S$ as above is an $R$-algebra, and that 
$$
E_\bullet: \quad \ldots \rTo E_{i+1} \rTo E_i \rTo \dots \rTo E_0
$$ 
as above
is the linear part of  $S\tensor_R F_\bullet$.  If the homology
$\H_{1}({E}_\bullet)$ is nonzero only in degrees $>1$, 
then ${E}_\bullet$ is linearly exact.

\noindent{\sl Proof.\ } We must show that if 
$\Tor^R_1(S,M)_1=0$, then $\Tor^R_i(S,M)_i=0$ for each $i\geq 1$;
a ``linear rigidity'' theorem for Tor.
The proof of Auslander-Buchsbaum
[1958] for the rigidity of Tor (reduction to the diagonal plus
the rigidity of the Koszul complex) may easily be adapted.\Box

In the theory above it actually suffices to suppose that $E_i$ is
generated in degree $i$ just for $i>n$. Thus we may try to apply the 
theory to the Eagon-Northcott complex and so we see that the complex 
is linearly exact iff the minors are independent. This leads to the
following:\smallskip

\noindent{\sl Problem:\ } Under what conditions are 
the $d\times d$ minors of an $e\times f$ matrix of linear forms linearly
independent?

\beginsection 3.  The complexes $E^m_\bullet(\mu)$

In this section we define the complexes whose linear exactness  plays a role in 
our analysis of the resolutions of ideals of points.  With appropriate
choices, these complexes extend
the map defined in Theorem 1.4 (for the moment only in the case $n=2$)
and thus admit a map to the dual of the free resolution of an ideal of 
points.

First we recall the notion of divided power.  Let $U$ be a 
finitely generated free module over some ring.
We write $D_l U$ for the $l\th$ divided power of $U$.  It is convenient 
to define $D_l U$ as the dual of the $l\th$ symmetric
power of the dual module, that is $D_l U = (\Sym_l(U^*))^*$.  What we shall use
about $D_l U$ is that it has a ``diagonal'' map $D_{l+1} U \rTo D_l U\tensor U$
which is the monomorphism
dual to the surjective multiplication map 
$\Sym_l (U^*)\tensor U^* \rTo \Sym_{l+1} (U^*)$. See for example
Eisenbud [1995, Appendix 2] for the usual definition.

Suppose again that
$S$ is a graded ring, with $S_0=k$ a field,
and that $U$ and $W$ are finite dimensional vector spaces over $k$.  Let
$\mu: W\tensor U\to S_1$ be a homomorphism.

For any integer $m$, and any integer $l\geq 0$ we define a free module 
$$
E^m_l(\mu):=\wedge^{l-m}W\tensor D_{l}U\tensor S(-l)
$$ 
and a map
$$
\delta^{m}_{l+1}(\mu):\ E^m_{l+1}(\mu) \rTo E^m_l(\mu),
$$
which is the composite of the tensor product of the diagonal maps
of the exterior and divided powers,
$$
\wedge^{l+1-m}W\tensor D_{l+1}U\tensor S(-l-1)\rTo 
\wedge^{l-m}W\tensor W\tensor D_{l}U\tensor U\tensor S(-l-1),
$$
and the map induced by $\mu$
$$
\wedge^{l-m}W\tensor W\tensor D_{l}U\tensor U\tensor S(-l-1)
\rTo 
\wedge^{l-m}W\tensor D_{l}U\tensor S(-l).
$$
These maps form complexes of free $S$-modules
$$
E^m_\bullet(\mu):\qquad
\ldots \rTo 
E^{m}_{l+1}(\mu)\rTo^{\delta^{m}_{l+1}(\mu)}
E^{m}_{l}(\mu)
\rTo
\ldots
\rTo 
E^{m}_{m_+}(\mu)
$$
where the term $E^{m}_{l}(\mu)$ is in position $l$,
and $m_+$, which denotes the positive part of $m$, is equal to $m$ if
$m\geq 0$ and to 0 if $m\leq 0$. The complex $E_\bullet(\mu)$ in
the introduction is $E^{-1}_\bullet(\mu)$ in the present notation.

The complexes $E^m_\bullet(\mu)$ are
linear in the sense above. If $u$ is the rank of $U$,
then $E^{-u}_\bullet(\mu)$ is precisely the linear part of the
Eagon-Northcott complex resolving the maximal minors of $\mu$,
whence the name $E$. As with the Eagon-Northcott
complex, these complexes may be built up inductively:

\proclaim Proposition 3.1.(Inductive Construction). With notation as above, 
suppose that $$0\rTo W^\prime \rTo W \rTo k\rTo 0$$ is an exact sequence, 
and let $\mu':W'\tensor U\rTo S_1$ denote the composition of $\mu$ with the
inclusion $W'\tensor U \rTo W\tensor U$.  There is
an exact sequence of complexes
$$
0 \rTo 
E^m_\bullet(\mu')\rTo E^m_\bullet(\mu) \rTo E^{m+1}_\bullet(\mu')
\rTo 0.
$$

\noindent{\sl Proof.\ } We use the exact sequence
$0\rTo \wedge^{l-m}W'\rTo \wedge^{l-m}W\rTo \wedge^{l-m-1}W'\rTo 0$.
The commutativity of the necessary diagrams follows by straightforward
computation.\Box

Using Proposition 3.1 we can show that the complexes $E_\bullet(\mu)$ satisfy the
hypothesis of the linear rigidity lemma above.

\proclaim Theorem 3.2. Let $k$ be a field, and let
$W$ and $U$ be finite dimensional vector spaces over $k$.
Let $S=k[W\tensor U]$ be the symmetric algebra, and let 
$\mu:\ W\tensor U \rTo S(1)$ be the identity map.
\item{a)} The complex $E^m_\bullet(\mu)$ is the linear part of a  minimal
free resolution. 
\item{b)} For any integers $m$
and $i>m_+$ the module $\H_i(E^m_\bullet(\mu))$ 
is nonzero at most in degrees $>i$.

\noindent{\sl Proof.\ }  An argument similar to that of 
Proposition 2.1 shows that parts {\sl a)} and {\sl b)} are equivalent.

To prove part {\sl b)}, we do induction on the rank $w$ of $W$.  If $w=1$ 
and $m<0$ there is nothing to prove. If $w=1$ and $m\geq 0$, then
exactness follows from the fact that the diagonal map
$D_{m+1}U\rTo D_{m}U\tensor U$ is a monomorphism.

Suppose now $w>1$. Let $W'$ be a codimension 1 subspace of $W$,
so that we have an
exact sequence $0\rTo W' \rTo W \rTo k \rTo 0$.
Using the long exact sequence in homology coming from
the inductive construction in Proposition 3.1, everything is
clear except the cases where $m\geq 0$ and $i=m+1$. In this case the
exact sequence of complexes has the form
$$
\diagram[tight,width=6em,height=2em]
{\bf Homological\ degree:}&&&{\bi m}+1&&{\bi m}\\
E^{m+1}_\bullet(\mu'):\qquad&\quad\qquad\ldots&\rTo&D_{m+1}U\tensor S(-m-1)&& \\
\uTo&&&\uTo&& \\
E^{m}_\bullet(\mu):\qquad&\quad\qquad\ldots&\rTo&W\tensor D_{m+1}U\tensor S(-m-1)&
                               \rTo&D_mU\tensor S(-m) \\
\uTo&&&\uTo&&\uEquals \\
E^{m}_\bullet(\mu'):\qquad&\quad\qquad\ldots&
                  \rTo&W'\tensor D_{m+1}U\tensor S(-m-1)&\rTo&D_mU\tensor S(-m) \\
\enddiagram
$$
and we must show that the connecting homomorphism
$$
c:\quad H_{m+1}(E^{m+1}_\bullet(\mu'))\rTo H_m(E^m_\bullet(\mu'))
$$
is a monomorphism in degree $m+1$, which is the lowest degree
present in $E^{m+1}_{m+1}(\mu') = D_{m+1}U\tensor S(-m-1)$.  

We have
$(H_{m+1}(E^{m+1}_\bullet(\mu')))_{m+1} = D_{m+1}(U)$. 
Let $f\in D_{m+1}(U)$, and write 
$\sum_i f_i\tensor f_i' \in U\tensor D_m(U)$
for the image of the diagonal map.  
If we write $W =\langle x\rangle\oplus W'$, and 
$\bar f$ for the class of $f$ in the homology of 
$E^{m+1}_\bullet(\mu')$, then we see by chasing the diagram that
$$
c(\bar f)\quad =\quad
\sum_i (x\tensor f_i)\tensor f_i'\quad \in\quad (W\tensor U)\tensor D_m(U).
$$
Since the differential of 
$E^{m}_\bullet(\mu')$ involves only $W'$, the homology module
$\H_m(E^{m}_\bullet(\mu'))$ surjects onto
$(\langle x\rangle\tensor U)\tensor D_m(U) = U\tensor D_m(U)$,
and we may recover 
$\sum_i f_i\tensor f_i'$ as the image of $c(\bar f)$.  Since the diagonal
map is a monomorphism on the divided powers, we are done.\Box

\medskip

A closer examination of the induction shows that the generic complexes
in Theorem 3.2 actually are resolutions if $m\leq -w+1$, but not 
otherwise; for example, if $m=0,\ w = \dim W = 2$,
and $\dim U := u > 2$, then the
complex $E^m_\bullet(\mu)$ has the form
$$
S^u(-1)\cong \wedge^2 W\tensor U\tensor S(-1)\rTo W\tensor S = S^2,
$$
and the free resolution of which this is the linear part is the
Buchsbaum-Rim complex
$$
\ldots \rTo W^*\tensor(\wedge^3 S^u)(-3) \rTo S^u(-1)\rTo S^2.
$$
The degree 2 relations $ W^*\tensor(\wedge^3 S^u)(-3) 
\rTo S^u(-1)$ are an expression
of Cramer's rule. See Eisenbud [1995, Appendix A2.6] for more
information.

In our setting the map $\mu$ has a geometric origin, and we may use
a technique similar to Green's Koszul Homology to check the condition
of linear exactness.  The following is the result of this section that
we shall use in the sequel:

\proclaim Corollary 3.3. Let $C$ be a projective scheme over a field $k$, 
and let
$H, L$ be Cartier divisors on $C$.  Suppose that $\O_C(H)$ is generated
by its global sections $W := \H^0(\O_C(H))$, 
and let $M_H$ be the vector bundle on $C$ which is
the kernel of the evaluation
map $d_0:\ \O_C\tensor W \rTo \O_C(H)$. Set $U:=\H^0(\O_C(L))$, and
$V:= \H^0(\O_C(H+L))$.  Let $S= \Sym\, V$ be the polynomial ring, and
let $\mu:\ W\tensor U \rTo V = S_1$ be the multiplication map.
The complex $E^{-1}_\bullet(\mu)$ is linearly exact if and
only if $H^0(\wedge^2 M_H\otimes \O_C(L))=0$.

\noindent{\sl Proof.\ } By the linear rigidity lemma it is enough
to check linear exactness at the first step; that is, we must
show that the induced map $\wedge^2 W\tensor U\rTo W\tensor V$
is a monomorphism. For this purpose
we use the Koszul complex built on the evaluation map $d_0$,
$$
\dots\rTo^{d_3}\wedge^3 W\otimes\O_C(-3H)\rTo^{d_2}\wedge^2 W\otimes\O_C(-2H)
\rTo^{d_1} W\otimes\O_C(-H)\rTo^{d_0}\O_C\rTo 0,
$$
tensored with $\O_C(2H+L)$.
Now $\ker d_i=\Im d_{i+1}=\wedge^{i+1}M_H(-(i+1)H)$, for all $i\ge 0$, so
the claim of the lemma follows by taking global sections in the  
short exact sequence
$$0\rTo\wedge^2 M_H\otimes\O_C(L)\rTo\wedge^2 W\otimes\O_C(L)\rTo^{d_1\tensor
\O_C(2H+L)} W\otimes \O_C(H+L).$$
\Box

\medskip
\noindent{\sl Remark.\/} We have made the restriction to projective schemes
only to ensure the finite dimensionality of the spaces involved.  This 
is actually unnecessary; the complexes $E^m_\bullet$ could have been developed
for infinite dimensional spaces.  We leave these things to the reader who can
find an application \dots

\beginsection \S 4. Subcomplexes of the resolution of $I_\Gamma$

We prove in this section the main result concerning resolutions of
points: For a
suitably chosen map $\mu$ the  complex $E^{-1}_\bullet(\mu)$ defined above  
is linearly exact, and its dual is a subcomplex 
of the back end of the minimal free resolution 
of the ideal of the points. 

\proclaim Theorem 4.1. Let $V$ be an $(r+1)$-dimensional vector space
over a field $k$ of characteristic 0, and let $\Gamma$ be a general set of
$\gamma$ points in $\P(V)$. Let $W:=V^\perp\subset \H^0( K_\Gamma(-1))$,
let 
$U:=(W^2)^\perp\subset \H^0(K_\Gamma^{-1}(2))$, and
let $\mu:\ W\tensor U \rTo V$ be the multiplication map.
Set 
$$
r:={{s+1\choose 2}+t},\ \ s\ge 2, \ \ 0 \leq t \leq s, 
{\rm \ \ and\ }\gamma := r+s+2.
$$
If $s$ is even suppose also ${\rm char}\, k=0$.
The complex $E_\bullet^{-1}(\mu)(r+2)$ is a direct summand of the
dual of the free resolution of the ideal of $\Gamma$.

{\noindent \sl Remarks.\ }The given number of points in $\P^r=\P(V)$ 
is actually the largest number for which the construction is interesting;  
for smaller numbers there is still a nontrivial complex but it is only 
sometimes  linearly exact. The restriction to characteristic 0 is most likely 
unnecessary, but is not important as our main Theorem 0.1 follows
in all characteristics from the characteristic 0 case. The restriction
comes only from the use of a theorem of Hartshorne and Gieseker on the
semistability of symmetric powers of semistable vector bundles
at the very end of the argument.
\medskip

{\noindent \sl Proof.\ } We shall  show that
the complex $E_\bullet^{-1}(\mu)$ is linearly exact.
Since $s\geq 2$ we have $\gamma\leq {r+2 \choose 2}$,
so $\Gamma$ imposes independent conditions on quadrics and
thus the homogeneous ideal $I_\Gamma$ is 3-regular.  
It follows that, with notation as in 
Theorem 1.4, $\omega_{\Gamma, V}$ is generated in degrees $\geq -1$.
The dual of the free resolution of $I_\Gamma$ is (the beginning of)
the minimal free resolution of $\omega_{\Gamma,V}(r+1)$.
We will thus deduce Theorem 4.1 from Theorem 1.4 and
Proposition 2.1, applied to
the complex $E_\bullet^{-1}(\mu)$ and the minimal free resolution of 
$\omega_{\Gamma,V}(-1)$.

Our strategy for proving linear exactness is as follows.
We wish to apply Corollary 3.3.  To do this we must find a scheme
$C$ such that $W$ may be interpreted as a space of sections generating
a line bundle $\O_C(H)$ and $U$ may be interpreted as the space of 
all sections of a line bundle $\L$. It is most convenient
to regard $\Gamma$ by its ``other'' embedding as the Gale transform
$\Gamma'$, since there $W$ is the space of sections of the line bundle
responsible for the embedding in $\P^s$, while $U$ may be identified
with $\H^1(\I_{\Gamma^\prime}(2))^*$. We cannot take 
$C=\Gamma^\prime$ itself, however,
because $U$ is not a complete linear series.  Thus we need some
higher-dimensional scheme on which $\Gamma'$ lies.  Since the general set
of points $\Gamma'$ does not 
(as far as we know) lie on any useful schemes of larger dimension,
we will make a degeneration, using the (obvious) openness of the
locus, in the space of maps $\mu: W\tensor U\rTo V$, where
$E^{-1}_\bullet(\mu)$ is linearly exact.
We shall degenerate $\Gamma'$ to a set of points,
lying on a convenient curve $C$.  In doing this, 
we must keep the dimensions of $W$ and $U$ constant
(since $V=W^\perp$, the constancy of its dimension is then automatic).

Since $\Gamma'$ is a general set of $\gamma > {s+2\choose 2}$ points in $\P^s$,
it lies on no quadrics, and this fact determines the dimension of $U$
as
$h^1(\I_{\Gamma^\prime}(2))$. We may thus degenerate $\Gamma'$
to a general subset of a curve $C$ in $\P^s$ that lies on no 
quadrics (which we will again call
$\Gamma^\prime$). 

In order to establish a simple relation between the cohomology
of $\I_{\Gamma^\prime}$ and bundles on the curve we will require 
$C$ to be nonspecial and quadratically normal.
Thus writing $H$ for the hyperplane class on $C$ and setting
$d:={\rm deg\ } H,\ g:={\rm genus\ }C$, we need
$$\displaylines{
s+1=h^0(\O_C(H))=d+1-g\cr
{{s+2}\choose 2}=h^0(\O_C(2H))=2d+1-g,\cr}
$$
which in turn yield $d={{s+1}\choose 2}$ and $g ={s\choose 2}$.

It is easy to compute that the curve defined in $\P^s$ by the 
vanishing of the $3\times 3$ minors of a general $3\times (s+1)$
matrix of linear forms $M$ has exactly the invariants required.
From the existence of this curve $C$, and the openness of the desired 
properties, we see that we may take $C$ to be a general curve of
genus ${s\choose 2}$, embedded by the complete linear series
associated to a general divisor $H$ of degree ${{s+1}\choose 2}$
in $\P^s$.
We will use this freedom to make further degenerations.

The binomial form of the genus formula suggests a plane curve of 
degree $s+1$, and it is amusing to note that the determinantal curve
just defined may be embedded in the plane by the line bundle that is
the cokernel of the restriction of $M$ to the curve; in this
planar embedding its equation is the determinant of the $(s+1)\times(s+1)$
matrix of linear forms in 3 variables which is adjoint to $M$. We shall
use this construction implicitly later in the proof.

If $\Gamma^\prime$ is a general divisor
of degree $\gamma$ on a curve $C$ as above, then 
we can write $\mu$ as a map coming from bundles on $C$ as follows:
Since $\O_C(H)$ is nonspecial and the curve $C$ is projectively normal,
the cohomology of the short exact sequences 
$$
0\rTo
\I_C(mH)
\rTo
\I_{\Gamma^\prime}(mH)
\rTo
\O_C(mH-\Gamma^\prime)
\rTo 0,
$$
together with Serre duality yield
$$\H^1(\I_{\Gamma^\prime}(mH))\iso
\H^1(\O_C(mH-\Gamma^\prime))\iso
(\H^0(\O_C(K_C+\Gamma^\prime-mH)))^\ast,\quad {\rm for\ all\ } m\ge 1.
$$
We now set  $L:=K_C+\Gamma^\prime-2H$, and we have
$U = \H^0(\O_C(L))$ as required. Since 
$\gamma$ is greater than the genus of $C$, we
may simply describe $L$ as the general divisor of degree
$2g-2+\gamma-2d = r-s\leq {s+1\choose 2}$. With
these identifications the pairing $\mu:\ W\tensor U \rTo V$ 
becomes the multiplication 
$$
\mu: \H^0(\O_C(H))\otimes \H^0(\O_C(L))
\rTo \H^0(\O_C(L+H)).
$$

By Corollary 3.3 it now suffices to prove for each $s\geq 2$ that
$\H^0(\wedge^2 M_H\otimes\O_C(L)) = 0$ where
\item{$\bullet$} $C$ is a general curve of genus $g :={s\choose 2}$.
\item{$\bullet$} $H$ is a general divisor on $C$ of degree $d:={s+1\choose 2}$.
\item{$\bullet$} $L$ is a general divisor on $C$ of degree $\leq l:= {s+1\choose 2}$.
\hfill\break

The standard method of proving such vanishing is by filtration and
stability (see Green-Lazarsfeld [1986], 
Ein-Lazarsfeld [1992]) but it does not yield a 
strong enough result, and the stability of $M_H$ 
would not be a strong enough condition,
so instead we shall use further degenerations:
Depending on the parity of $s$, we reduce
to the case where $M_H$ is a direct sum of line
bundles ($s$ odd), or a direct sum of rank 2 vector bundles ($s$ even).
The desired vanishing is an open condition on 
the triples $(C,H,L)$ in any flat family for which
the dimension of $H^0(\O_C(H))$ is constant. 
Suppose that we can find, for each $s$,
a smooth curve $C_0$ of genus $g$, 
a nonspecial divisor $H$ of degree $d$, and some divisor $L$ of
degree $l$ on $C$
such that $\H^0(\wedge^2 M_H\otimes\O_C(L)) = 0$. 
Over a versal deformation of
$C_0$ we may form the space of triples $(C, H, L)$, where 
$H$, $L$ are
divisors of the given degrees. The base space of this 
versal deformation
maps to the moduli space of curves of genus 
$g$ and covers an open set
therein. Thus the general curve $C$, with general 
divisors $H$ and $L$
will have the properties required.
\medskip

Assume now that $s$ is odd. We 
let  $C_0$ be a curve of type
$((s+1)/2,s+1)$ on the smooth quadric $Q\subset\P^3$.
Let 
$\N$ be the restriction of $\O_Q(0,1)$ to $C_0$.
Thus $\N$ is a line bundle of degree ${(s+1)}/2$ generated
by global sections, and $\N^{\oplus s}$ is a globally generated
vector bundle of rank $s$ and degree ${s+1\choose 2}$ on $C_0$. 
Let $W^\ast$ be a general $(s+1)$-dimensional subspace of sections 
$W^\ast\subset H^0(\N^{\oplus s})$.  It is easy to see that
$W^\ast$ generates $\N^{\oplus s}$, and so we define a
line bundle $\O_{C_0}(H)$ as the dual of the
kernel of the natural evaluation
$$
0\rTo\O_{C_0}(-H)\rTo W^\ast\otimes\O_{C_0}\rTo \N^{\oplus s}\rTo 0.
$$

With these choices $\O_{C_0}(H)$ is a globally generated 
line bundle of the desired degree ${s+1\choose 2}$, and $W$
maps naturally to $\H^0(\O_{C_0}(H))$. In fact
$\O_{C_0}(H)=\det(\N^{\oplus s})=\N^{\otimes s}=\O_{C_0}(0,s)$. Furthermore,
since $h^1(\O_Q(0,s))=0$ and $h^2(\O_Q(-{{(s+1)}\over 2},-1))=0$,
taking cohomology of the exact sequence
$$0\rTo\O_Q(-{{(s+1)}\over 2},-1)\rTo\O_Q(0,s)\rTo\O_{C_0}(H)\rTo 0$$
we get $h^1(\O_{C_0}(H))=0$, that is
$\O_{C_0}(H)$ is non-special.
Thus $h^0(\O_{C_0}(H))=\chi(\O_{C_0}(H))=s+1$ and 
since $h^0(\N^\ast)=0$ for degree reasons, we get that
$W=H^0(\O_{C_0}(H))$, whence
$M_H\iso (\N^{\oplus s})^\ast$.
To show that $\H^0(\wedge^2 M_H\tensor\O_C(L))$, it is enough
to show that $\H^0(\N^{-2}\otimes\O_C(L))=0$. This is obviously true
for a general $L$ with $\deg L=r-s$, since
$\deg(\O_C(L)\otimes \N^{-2})=r-s-2\le {s\choose 2}-1=
g({C_0})-1$ by our initial hypothesis.

Consider the versal deformation of the curve
${C_0}$ and over it the space of triples $(C, H, L)$ as above, where
$H$ is a divisor of degree $d$ and $L$ is a divisor of degree
$r-s$.
The locus for which $\O_C(H)$ defines
an arithmetically normal embedding in $\P^s$ is 
open and, as we have seen, non-empty. Furthermore, the vanishing of 
$\H^0(\wedge^2 M_H\otimes \O_C(L))=0$ is an open condition on the
collection of triples.  Since the vanishing condition
is satisfied on $C_0$, the same follows for the general curve.

\bigskip

Finally, consider the case where $s$ is even. In order to produce
a nonspecial divisor $H$ with the desired properties in this case,
we will degenerate further, letting $H$ become special. Thus we
must work with incomplete linear series.

Given a divisor $H$ on a curve $C$ and a space of global sections
$W$ that generates  $\O_C(H)$, we define
$M_{W,H}$ to be the kernel of the natural
evaluation map:
$$
M_{W,H}:=\ker( W\otimes\O_C\rOnto\O_C(H)).
$$
It now suffices, for each even $s$, to find:
A curve $C$ of genus $g$, a divisor $H$ of degree
$d$ on $C$, and a space of sections $W$ of dimension $s+1$ of $\O_C(H)$
such that 
\item{$\bullet$} $\H^0(\wedge^2 M_{W,H}\otimes\O_C(L)) = 0$ 
for the general divisor $L$ on $C$ of degree $r-s$, and
\item{$\bullet$} The triple $(C, H, W)$ is a flat limit of  triples
for which $H$ is nonspecial (equivalently, where $W=\H^0(\O_C(H))$).

A candidate is constructed for us by the following result:

\proclaim Proposition 4.2.  Let $s\geq 2$ be even. For
any sufficiently general plane curve $C_0$ of degree
$s+1$ there  exists a flat irreducible
family of degree smooth plane curves $C_t$, with special fiber
$C_0$ and general fiber $C_\eta$, a family of line bundles $H_t$
of degree ${{s+1}\choose 2}$, and a family of spaces 
$W_t\subset \H^0(\O_{C_t}(H_t))$ such that:
\item{a)}  $W_t$ generates $\O_{C_t}(H_t)$,
\item {b)} $H_\eta$ is nonspecial, and
\item {c)} $M_{H_0}:= \ker (W_0\tensor \O_{C_0}\rTo \O_{C_0}(H_0))$
is the direct sum of $s/2$ copies of the rank 2 vector bundle $M$
which is the kernel of the evaluation map 
$\H^0(\O_{C_0}(1))\tensor\O_{C_0}\rOnto \O_{C_0}(1)$, where
$\O_{C_0}(1)$ induces the planar embedding.

\noindent{\sl Proof.\/} We shall construct the family of curves $C_t$
and the family of divisors $H_t$ by constructing the family of
vector bundles  
$\E_t := \ker (W_t\tensor \O_{C_t}\rTo \O_{C_t}(H_t))$.
On the generic fiber, we use the following (old) observation:  If 
$B:\ \O_{\P^2}^{s+1}(-1)\rTo \O_{\P^2}^{s+1}$ is an
$(s+1)\times(s+1)$ matrix of linear forms in 3 variables such that
$f := \det B \neq 0$, and such that the ideal $I_s(B)$ of 
$s\times s$-minors of $B$ contains a power of the irrelevant ideal, then
${\cal H} := \coker B$ is a line bundle on the degree $s+1$ curve $\{f=0\}$ with 
$\h^0({\cal H}) = s+1,\ \h^1({\cal H}) = 0$ (and thus $\deg {\cal H} = 
{s+1\choose 2}$ by Riemann-Roch). The vector bundle 
$M_{\cal H}$ is of course the image of $B$.

For the special fiber, we proceed differently.
Recall that  $\widetilde M := \Omega^1_{\P^2}(1)$ is the 
image of the middle Koszul map 
$$
\rho:\wedge^2 \H^0(\O_{\P^2}(1))\otimes
\O_{\P^2}(-1)\rTo \H^0(\O_{\P^2}(1))\otimes\O_{\P^2},
$$
induced by the $3\times 3$ generic skew-symmetric 
matrix  over the ring $S = \Sym(\H^0(\O_{\P^2}(1)))\iso k[x,y,z]$. 
If $C_0$ is any plane curve, then the bundle $M$ defined
in the Proposition is simply $\widetilde M |_{C_0}$.  We wish to 
define a general matrix of linear forms whose image is $\widetilde M$.
For this (and for later purposes) the idea of a ``generalized 
submatrix'' of a matrix will be useful:  by a generalized $p\times q$
submatrix of a matrix $C$ we mean simply a composition 
$PCQ$, where $P$ and $Q$ are scalar matrices, $P$ has $p$ rows,
and $Q$ has $q$ columns. Generalized rows or columns of $C$ are 
generalized submatrices with one row or one column, respectively.

Now let $A$ be a (sufficiently general)
generalized $(s+1)\times(s+1)$-submatrix 
of a  $(3s/2)\times (3s/2)$-matrix inducing
$\rho^{\oplus {s\over 2}}$. Notice that
$\det A=0$ since the module ${\Omega^1_{\P^2}(1)}^{\oplus {s\over 2}}$ 
has only rank $s$. Let
$B$ be a general $(s+1)\times(s+1)$-matrix with linear entries, and set 
$A_t := A+tB$. For $t\neq 0$, we set $f_t := \det A_t$, and for
$t=0$ we take $f_t$ to be the ``limit''
$$
f_0:=\lim_{t\rightarrow 0}{{\det(A+t\cdot B)}\over t}=
\sum_{i,j=1}^{s+1} b_{ij}\; |A_{ij}|,
$$
where $b_{ij}$ are the entries of $B$, while 
$|A_{ij}|$ denotes the (signed) minor of $A$ obtained by deleting row
$i$ and column $j$.

\proclaim Proposition 4.3. Set $m=s/2$. For each generalized
$(s+1)\times (s+1)$ submatrix $A$ of $\rho^{\oplus m}$, 
the ideal of $s\times s$ minors of $A$ may be written in the form
$I_s(A) = I_{m-1}(K_1)\cdot I_{m-1}(K_2) \cdot (x,y,z)^2$ for matrices of
linear forms $K_1$, $K_2$ of sizes $m\times (m-1)$ and $(m-1)\times m$,
respectively.  Each pair of matrices $K_1$ and $K_2$ with linear entries
arises for some generalized submatrix $A$. 
Thus, for a general choice of $A$, the ideal $I_s(A)$ is
a nonsaturated ideal of a reduced set of points, and the equation 
$\{f=0\}$ defines a general plane curve $C$ of degree $(s+1)$.

{\sl Proof of Proposition 4.3.\ } 
To compute $I_s(A)$ we make use of a special case of
the structure theorem for finite free resolutions of
Buchsbaum-Eisenbud [1974]. Consider the resolution
obtained by taking the 
direct sum of $m$ copies
of the Koszul complex in 3 variables:
$$
0
\rTo
S^{m}(-2)
\rTo^{{\kappa^\ast(-1)}^{\oplus m}}
S^{3m}(-1)
\rTo^{\rho^{\oplus m}} 
S^{3m}\rTo^{\kappa^{\oplus m}}
S^{m}(1)
\rTo 
0;
$$ 
to simplify notation, write it as
$$0\rTo F_3\rTo^{f_3} F_2\rTo^{f_2} F_1\rTo^{f_1} F_0.$$
Let  $r_i:=\rank(f_i)$, so that in particular 
$\rank(F_i)=r_i+r_{i+1}$, for all
$1\le i\le 3$.

The structure theorem 
asserts the commutativity (up to a sign) of the 
diagram:
$$
\diagram[tight,width=4em,height=3em]
\wedge^{r_3} F^\ast_2&=\wedge^{r_2} F_2&& 
\rTo^{\wedge^{r_2}f_2}&& \wedge^{r_2}F_{1}=&\wedge^{r_1}F_{1}^\ast\\
&\rdTo_{\wedge^{r_3}f_3^\ast}&&&&\ruTo_{\wedge^{r_1}f_1^\ast}\\
&&\wedge^{r_3}F_3^\ast&\cong R\cong &\wedge^{r_1}F_0^\ast&&\\
\enddiagram
\eqno{(\ast)}
$$
In other words, any minor of order $r_2$ of $f_2$ may be 
expressed as the product of complementary minors of
orders $r_1$ and $r_3$ of $f_1$ and $f_3$, respectively.

The choice of the matrix $A$ involves the choice of $s+1 = 2m+1$
generalized rows and
columns  of the matrix defining 
$\rho^{\oplus m}$, hence
the choice of $m-1$ complementary columns of $\kappa$ and $m-1$ 
complementary rows of $\kappa^\ast$, respectively.  
We denote by $K_1$ and $K_2$ 
the $m\times (m-1)$ and $(m-1)\times m$-submatrices 
of $\kappa$ and $\kappa^\ast(-1)$
distinguished in this way. 

Because of the structure
of  $\kappa$, any $m\times(m-1)$ matrix
with linear entries in $S$ can be obtained as $K_1$
through the appropriate choice
of $(m-1)$ generalized columns of $\kappa$, and similarly for
$K_2$. In particular, $K_i$ may be chosen to make
$I_{m-1}(K_i)$ be the ideal of any sufficiently general set
of $m\choose 2$ points in the plane. Diagram $(\ast)$
expresses the $(2m)\times (2m)$-minors $A_{ij}$ of $A$ as the
products of $m\times m$ minors
of $\kappa$ and $\kappa^*$ that contain $K_1$ and $K_2$, respectively. 
An $m\times m$-minor of $\kappa$ containing $K_1$
is a linear combination of the $(m-1)\times (m-1)$
minors of $K_1$, with coefficients the elements of an arbitrary generalized
column of $\kappa$.  Again because of the structure of $\kappa$
this column can be taken to be an arbitrary column of linear forms
in  $S$.  Thus the ideal of $m\times m$-minors of $\kappa$ containing
$K_1$ is $I_{m-1}(K_1)\cdot (x,y,z)$. As similar remarks hold for
$K_2$, we have proven the first part of the Proposition.

If the choice of the generalized submatrix $A$ is general,
then $K_1$ and $K_2$ will be general matrices of linear forms,
hence their ideals of minors will be reduced ideals of distinct
general sets of points in the plane, and the ideal
$I_s(A) = I_{m-1}(K_1)\cdot I_{m-1}(K_2) \cdot (x,y,z)^2$
will be a nonsaturated ideal of the union of these two sets of 
points, as claimed.

Varying the
matrix $B$ we obtain for $f$ any form of degree $(s+1)$
in  $I_s(A)(x,y,z) = I_{m-1}(K_1)\cdot I_{m-1}(K_2) \cdot (x,y,z)^3$. 
Since ${s+3\choose 2} > 2{m\choose 2}$ the general curve
$C_0$ of degree $s+1$  through two general sets of $m\choose 2$
points in the plane is a general plane curve of degree $s+1$,
concluding the argument.
\Box

\noindent
{\sl Completion of the proof of Proposition 4.2.} 
Let $C_t$ be defined  by the equation $\{f_t=0\}$, and
let $\E_t$  be the image of the restriction 
to  $C_t$ of the
morphism induced by the matrix $A+t\cdot B$. Let
${\cal H}_t^{-1}$ be the kernel of the  restriction of the matrix
$(A+Bt)^*$ to $C_t$; we write the dual in the form
$\O_{C_t}(H_t) = {\cal H}_t$, for a family of divisors $H_t$
(defined, for example, by the family of sections that are the images
of the first basis vector of the target free module of $A+Bt$).
Part {\sl a)} of the Proposition now follows from the definitions;
part {\sl b)} follows from the remark at the beginning of the proof;
and part {\sl c)} follows from the form of the matrix $A = A+0\cdot B$.
\Box

\noindent
{\sl Continuation of the proof of Theorem 4.1.} We adopt the notation of
Proposition 4.2, but for simplicity we now set $C = C_0$.
By Proposition 4.2,
$C$ may be chosen to be a general plane curve of degree $s+1$. 
It suffices
to show that
$\H^0(\wedge^2 (M^{\oplus s/2})\otimes\O_C(L))=0$,
where $L$ is a general divisor of degree $r-s$, and  
for this it is enough to show that both 
$\H^0(\wedge^2 M\otimes\O_C(L))=0$, and 
$\H^0(M\otimes M\otimes\O_C(L))=0$. 

Two remarks will make the plausibility of this conclusion clear.
First, $r\geq {s+1\choose 2}$ so $\deg(L) = r-s\geq {s\choose 2}$,
and $g = {s\choose 2}$ is the genus of $C$. Thus $\O_C(L)$ is a general 
line bundle in the Picard  variety of $C$. Second, $\deg(M) = -(s+1)$, so
$\chi(\wedge^2 M\otimes\O_C(L))$ and
$\chi(M\otimes M\otimes\O_C(L))$ are both $\leq 0$.  Thus
each of the desired vanishings has the form: $\H^0(F\otimes \O_C(L'))=0$
with $F$ a vector bundle on $C$ with $\chi(F)=0$, and $L'$  a general
divisor of degree $\leq 0$. This condition obviously implies
that the bundle $F$ must be semistable, and indeed Raynaud [1982]
shows that the condition is equivalent to semistability when 
$\rank F\leq 2$, and also when $\rank F = 3$ on a general curve of a given
genus. In fact his argument proves a little more:

\proclaim Theorem 4.4 (Raynaud). Let $C$ be a general plane
curve  of any degree $\geq 3$.
A vector bundle $F$ of rank $\leq 3$ on $C$ with $\chi(F) = 0$ is
semistable iff $\H^0(F\otimes \O_C(L')) = 0$ for the general line bundle 
$\O_C(L')$ of degree 0 on $C$.

\noindent{\sl Discussion of Theorem 4.4.\/} Raynaud [1982, \S 2] 
enunciates the result for general curves (not planar).  However, his proof
shows that if we replace ``vector bundle'' by ``torsion-free sheaf", then 
the truth of the Theorem for $C$ defines an open set in the 
moduli of stable curves. Furthermore, his proof shows that this open 
set includes every irreducible rational curve of arithmetic genus
$g$, having exactly $g$ ordinary nodes. Since the the general map from $\P^1$
into the plane has as image a curve with only ordinary nodes as singularities,
these facts imply that the Theorem holds for a general plane curve.\Box

The first of the necessary vanishings is immediate from the
remarks above: Since a general line bundle of degree $\leq g-1$ 
has no sections $\H^0(\wedge^2 M\otimes\O_C(L))=\H^0(\O_C(L-N))=0$, 
where $N$ is the divisor of the intersection of $C$ with a line.

For the second vanishing, from the exact sequence
$$
0\rTo \wedge^2 M\rTo  M\otimes M \rTo \Sym_2 M\rTo 0,
$$
together with the first vanishing result above, it suffices to show that
$\H^0(\Sym_2(M)\otimes\O_C(L))=0$, and this puts us in the case of a 
bundle $F = \Sym_2(M)\otimes\O_C(L))$ of rank 3. An easy degree
computation shows that $\chi(F)\leq 0$. We may now invoke
Theorem 4.4 to conclude the argument if we can show that 
$F$ is semistable, and by Hartshorne [1971], Gieseker [1979] 
it suffices to show that $M$ itself is semistable. 
As $M$ differs from the restriction
to $C$ of the tangent bundle of the projective plane only by
twisting by a line bundle, the following elementary result completes the
proof of Theorem 4.1:

\proclaim Proposition 4.5.  Let $T=T_{\P^2}$ be the tangent bundle of the
projective plane. If $C$ is a smooth
plane curve of degree $m\geq 3$, then $T|_C$ is stable.

\noindent{\sl Proof.\/} Suppose $Q$ is a line bundle  
quotient of $T|_C$. Since $\deg(T|_C)=3m$, we must 
show that $\deg(Q)>3m/2$. But $T(-1)$
is globally generated, so either $Q(-1) = \O_C$, or 
$Q(-1)$ defines a base point free linear series. Let
$e:=\deg(Q)-m$ be the degree of $Q(-1)$; we must show that 
$e>m/2$. 

Assume first that $Q(-1)=\O_C$. Restricting the presentation of 
$T$ to $C$ we obtain maps
$$ 
\O_C(-1) \rTo \O_C^3 \rOnto T|_C(-1) \rOnto Q(-1)= \O_C
$$ 
with composition 0, where the last two maps are surjective. Since the
restrictions of linear forms on $\P^2$ are still linearly independent on
$C$, this is a contradiction. 

Thus we may suppose that $Q(-1)$ defines a base point free linear series
of degree $e$. It follows at once that $e\geq m-1>m/2$. For the
reader's convenience we give the elementary proof:
Since the genus of $C$ is positive, we must have $e\geq 2$.  If
$m=3$, then this is the desired result. On the other hand, if $m\geq 4$,
then $m-2 < {m-1\choose 2}$, the genus of $C$, so if $e\leq m-2$
then $Q(-1)$ is special. In other words
the points in a divisor $D$ in the linear series represented by $Q$
impose dependent conditions on the canonical series
$\O_C(K_C)=\O_C(m-3)$. But any finite scheme of length $\leq m-2$ 
in the plane imposes independent conditions on forms of degree 
$m-3$. \Box

\beginsection \S 5. Numerology of resolutions and failure of 
the Minimal Rank Conjecture.

In this section we derive first a lower bound for the graded betti
numbers of the homogeneous coordinate ring of a scheme
$\Gamma$ of $\gamma$ general points in $\P^r$, and
then prove  Theorem 0.1, providing
counterexamples to the Minimal Resolution Conjecture. When
the lower bound is achieved, we shall say that 
$\Gamma$ has {\it expected betti numbers}.  It is 
well-known how to do the computation (it appears 
explicitly in the Queen's University thesis of Anna
Lorenzini as well as in Lorenzini [1987], [1993]), but 
for the reader's convenience,
and because we need details in a certain special case,
we spell it out.  Since all we 
use about $\Gamma$ is its Hilbert function, the same 
computation would work for any subscheme finite over $k$
imposing ``as many conditions as possible'' on 
forms of each degree; we call such a subscheme
{\it sufficiently general}. We shall use the following 
elementary facts:

\item {a)} If $S_\Gamma$ is the homogeneous
coordinate ring of a finite scheme $\Gamma$ of points
in $\P^r$, and $S$ is the homogeneous coordinate ring
of $\P^r$, and $x$ is a linear form not vanishing
on any point in the support of $\Gamma$, then
the graded betti numbers of $S_\Gamma$ as an $S$-module
are the same as the graded betti numbers of 
$S_\Gamma/xS_\Gamma$ as an $R := S/x$-module. 
\item{b)} Write 
$R=k[x_1,\ldots,x_r]$, set $m = (x_1,\dots,x_r)$,
and let $d$ be the largest integer such that
$\gamma \geq {r+d-1 \choose d-1}$. In other words,
assuming $\Gamma$ is sufficiently general,
$d$ is the smallest degree of a form contained
in the homogeneous ideal $I_\Gamma$.
We may write $S_\Gamma/xS_\Gamma$ in the form 
$R/I$ where $ m^{d+1} \subset I \subseteq m^{d}$.
\item{c)} The module $R/I$ may be obtained as a 
``jump deformation'' of the module
$R/m^{d} \oplus m^{d}/I$.  Thus the resolution of 
$R/I$ has all graded betti numbers $\geq$ the betti
numbers in the resolution of $R/m^{d} \oplus m^{d}/I$.

\smallskip
\noindent We get the desired estimates by 
putting these things together with a
knowledge of the free resolution of 
$R/m^d$, which may for example be described as an
Eagon-Northcott complex (see Eisenbud [1995]).  
We write $\tilde \beta_{i,j}$ for the
expected dimensions of the Koszul homology,
and $\left\{ n \right\}_+$ for 
$\max(n,0)$. We set
${n\choose k}=0$ for $k>n$.

\proclaim Proposition 5.1.  Let $\Gamma$ be a finite
sufficiently general subscheme of $\P^r$ having 
degree $\gamma$ with 
$$
{r+d-1 \choose d-1} \leq \gamma < {r+d \choose d},
$$
and set $a := \gamma - {r+d-1 \choose d-1} \geq 0$.
The Koszul homology dimensions
$$
\beta_{i,j}(I) =  \dim_k(\Tor^S_i(I_\Gamma, k)_{j})
$$
in the ``interesting'' range $0\leq i\leq r-1$ satisfy:
\item{a)}
$
\beta_{i,j} = 0 \ {\rm unless\ }j=i+d{\rm\ or\ }j=i+d+1;
$
\item{b)} 
$$\eqalign{
{d+i-1\choose i}{r+d-1\choose d+i} 
\geq \beta_{i,i+d}&
\geq \tilde\beta_{i,i+d} = \left\{
{d+i-1\choose i}{r+d-1\choose d+i}-a{r\choose i}
\right\}_+ \cr
a{r\choose i+1}
\geq \beta_{i,i+d+1}&
\geq 
\tilde\beta_{i,i+d+1} = \left\{
a{r\choose i+1}- {d+i\choose i+1}{r+d-1\choose d+i+1} 
\right\}_+
}
$$
\noindent\hfill\Box

For simplicity, and because it is the case we shall
use, we now specialize to the case where
$d=2$, so that the ideal of $\Gamma$ is generated by
quadrics and cubics. 

\proclaim Corollary 5.2.  Let $\Gamma$ be a finite
sufficiently general subscheme of $\P^r$ having 
degree $\gamma$ with $ r+1\leq \gamma < {r+2\choose 2}$.  The
expected
dimensions of the Koszul homology  of $\Gamma$ are: 
$$
\eqalign{
\tilde \beta_{i,i+2}=&
\left\{ (i+1){r+2\choose i+2}-\gamma{r\choose i}\right\}_+,\cr
\tilde \beta_{i,i+3}=&\left\{ \gamma{r\choose i+1}-(i+2)
{r+2\choose i+3}\right\}_+,
\quad i=\overline{0,r-1}\cr}
$$
In particular 
$
\tilde\beta_{i,i+2}\neq 0
$
iff
$i<{(r+2)(r+1)\over \gamma}-2$. Furthermore
$
\tilde\beta_{i,i+3}\neq 0
$
iff 
$i\geq {(r+2)(r+1)\over \gamma}-3$.

\noindent{\sl Proof.\ } Arithmetic, starting from the 
previous result. \Box

Thus the ``expected'' shape
minimal free resolution of $\I_\Gamma$ is

$$
\vbox{\offinterlineskip 
\halign{\strut\hfil# \ \vrule\quad&# \ &# \ &# \ &# \ &# \ &# \ 
&# \ &# \ &# \ &# \ &# \ &# \ &# \ &# \ 
\cr
degree&\cr
\noalign {\hrule}
0&1&--&--&$\ldots$&--&--&$\ldots$&--&--&--\cr
1&--&*&*&$\ldots$&*&?&$\ldots$&--&--&--\cr
2& & & &$\ldots$&?&*&$\ldots$&*&*&*\cr
\noalign{\smallskip}
\omit&\omit&\multispan{4}{\upbracefill\ }&&&&&&\cr
\noalign{\medskip}
\omit&\multispan{6}{\hfill{$\scriptstyle\bigl[{{(r+1)(r+2)}\over\gamma}
\bigr]-3$}\hfill}&&&&&\cr
\noalign{\smallskip}
}}
$$
(where not both of the ``?''s in the above display are non-zero !)

\medskip\noindent

As a an easy corollary of Theorem 4.1 on linear exactness we obtain now the 
result announced in the introduction.

\noindent{\sl Proof of Theorem 0.1.\ } 
 For $r$ and $s$ in the given range the
complex $E^{-1}_\bullet(\mu)$ defined at the beginning of section \S 4
is linearly exact. Moreover, the twisted complex $E^{-1}_\bullet(\mu)(r+2)$
maps monomorphically onto a direct summand  of the dual of the
minimal free resolution of $I_\Gamma$. On the other hand, Corollary 5.2
gives as expected graded betti number 
$\tilde\beta_{(r-s-1),(r-s+2)}$ for $I_\Gamma$ 
$$\tilde\beta_{(r-s-1),(r-s+2)}=
\{{{(2k+4-s^2+s)}\over{(s^2-s+2k+4)}}\cdot
{{s+1\choose 2}+k\choose {s\choose 2}+k}\}_+,$$
whereas the last (i.e., the $s$-th) syzygy 
module in the complex $E^{-1}_\bullet(\mu)$ has rank
$$\rank E^{-1}_{s-1}(\mu)={s+k\choose k}.$$
The theorem follows since
$$s^2-s+2k+4\ge s^2-s>0\quad{\rm and}\quad
2k+4-s^2+s\le 3s+4-s^2\le 0,\quad {\rm for\ all\ } s\ge 4,\ 0\le k\le s,$$
while
$${{(2k-2)}\over {(2k+10)}}\cdot {k+6\choose k+3}<{k+3\choose 3}$$
only for $r=k+6\in\{6,7,8\}$. \Box

\bigskip
{
\centerline {\bf References}
\baselineskip=12pt
\parindent=0pt
\frenchspacing
\medskip 

\item{} M.~Auslander, D.~Buchsbaum: Codimension  and multiplicity,
{\sl Annals of Math.} {\bf 68} (1958), 625--657.
\medskip 

\item{} E.~Ballico, A.V.~Geramita: The minimal free resolution of the ideal 
of $s$ general points in $\P^3$, 
{\sl Proceedings of the 1984 Vancouver conference in algebraic
geometry}, 1--10, CMS Conf. Proc., 6,  Amer. Math. Soc., 
Providence, R.I., 1986.
\medskip

\item{} D.~Bayer, M.~Stillman:
Macaulay: A system for computation in
        algebraic geometry and commutative algebra
Source and object code available for Unix and Macintosh
        computers. Contact the authors, or download from 
        {\bf math.harvard.edu} via anonymous ftp.
\medskip

\item{} S.~Beck, M.~Kreuzer: How to compute the canonical module of a set
of points, {\sl Proceedings of the MEGA conference, Santander 1994}
(to appear).
\medskip

\item{} M.~Boij:
Artin level algebras, Doctoral Dissertation, Stockholm 1994.
\medskip

\item{} D.~Buchsbaum, D.~Eisenbud:
Some structure theorems for finite free resolutions, {\sl 
Advances in Math.} {\bf 12}, (1974), 84--139.
\medskip

\item{} G.~Castelnuovo:
Su certi gruppi associati di punti,
{\sl Ren. di Circ. Matem. Palermo} {\bf 3},
(1889), 179--192.
\medskip

\item{} M.P.~Cavaliere, M.E.~Rossi, G.~Valla: On the
resolution of points in generic position, {\sl Comm.
Algebra} {\bf 19}, (1991), no. 4, 1083--1097.
\medskip

\item{} M.P.~Cavaliere, M.E.~Rossi, G.~Valla: 
The Green-Lazarsfeld conjecture for
$n+4$ points in $\P^n$,
{\sl Rend. Sem. Mat. Univ. Politec. Torino} {\bf 49},
(1991), 175--195. 
\medskip

\item{} A.B.~Coble: Point sets and allied Cremona groups I, II, III,
{\sl Trans. Amer. Math. Soc.} {\bf 16}, (1915), 155--198, 
{\sl Trans. Amer. Math. Soc.} {\bf 17}, (1916), 345--385,
and {\sl Trans. Amer. Math. Soc.} {\bf 18}, (1917), 331-372.
\medskip

\item{} A.B.~Coble: Associated sets of points,
{\sl Trans. Amer. Math. Soc.} {\bf 24}, (1922), 1--20.
\medskip

\item{} I.~Dolgachev, D.~Ortland: Points sets in
projective spaces and theta functions, {\sl Ast\'erisque}
{\bf 165}, (1988).
\medskip

\item{} L.~Ein, R.~Lazarsfeld:
Stability and restrictions of Picard bundles, with
an application to the normal bundles of elliptic curves, in
``Complex projective geometry (Trieste, 1989/Bergen, 1989)'', 149--156, 
{\sl London Math. Soc. Lecture Note Ser.}, {\bf 179}, Cambridge Univ.
Press, Cambridge, 1992.
\medskip

\item{} D.~Eisenbud:
{\sl Commutative Algebra with a View Toward Algebraic Geometry},
Springer, New York, 1995.
\medskip

\item{} D.~Eisenbud: Linear sections of determinantal varieties,
{\sl Amer. J. Math.} {\bf 110}, (1988) 541--575.
\medskip

\item{}D.~Eisenbud, S.~Popescu:  
The projective geometry of the Gale transform.
Preprint, 1997.
\medskip

\item{} F.~Gaeta: Sur la distribution des degr\'es des formes
appartenant \`a la matrice de l'id\'eal homog\`ene attach\'e
\`a un groupe de $N$ points g\'en\'eriques du plan,
{\sl C.~R.~Acad.~Sci.~Paris} {\bf 233} (1951) 912--913.
\medskip

\item{} F.~Gaeta: A fully explicit resolution of the ideal
defining $N$ generic points in the plane. Preprint, 1995.
\medskip

\item{}D.~Gale: Neighboring vertices on a convex polyhedron,
in ``Linear Inequalities and Related Systems'' (H.W.~Kuhn
and A.W.~Tucker, eds.), {\sl Annals of Math. Studies} {\bf 38},
255--263, Princeton Univ. Press, 1956.
\medskip

\item{} A.V.~Geramita, A.M.~Lorenzini:
The Cohen-Macaulay type of $n+3$ points in $\P^n$,
in ``The Curves Seminar at Queen's'', Vol. VI, Exp. No. F, 
{\sl Queen's Papers in Pure and Appl. Math.} {\bf 83}, (1989).
\medskip

\item{} D.~Gieseker:
On a theorem of Bogomolov on Chern classes of stable bundles. 
{\sl Amer. J. Math.}  {\bf 101}, (1979), no. 1, 77--85.
\medskip

\item{} V.D.~Goppa: Codes and Information, 
{\sl Russian Math. Surveys} {\bf 39}, (1984) 87--141.
\medskip

\item{} M.~Green, R.~Lazarsfeld: A simple proof of Petri's theorem
on canonical curves, in {\sl ``Geometry Today''}, Prog. in Math. Series,
Birkh\"auser (1986).
\medskip
 
\item{} M.~Green, R.~Lazarsfeld:
Some results on the syzygies of finite sets and algebraic curves,
{\sl Compositio Math.} {\bf 67}, (1988), no. 3, 301--314.
\medskip 

\item{} R.~Hartshorne:
Ample vector bundles on curves, 
{\sl Nagoya Math. J.} {\bf 43}, (1971), 73--89. 
\medskip

\item{} L.O.~Hesse: De octo punctis intersectionis trium
superficierum secundi ordinis,
Dissertatio, (1840), Regiomonti;
De curvis et superficiebus secundi ordinis,
{\sl J. reine und angew. Math.} {\bf 20}, (1840), 285--308.
\medskip

\item{} A.~Hirshowitz, C.~Simpson: La r\'esolution minimale
de l'arrangement d'un grand nombre de points dans $\P^n$,
to appear in {\sl Inventiones Math.}
\medskip

\item{}M.~Kreuzer: On the canonical module of $0$-dimensional scheme,
{\sl Canadian J. of Math.} {\bf 46} (1994), no. 2, 357--379.
\medskip

\item{}F.~Lauze: Rang maximal pour $T_{\P^n}$,
preprint alg-geom/9506010, (1995).
\medskip

\item{}F.~Lauze: preprint (in preparation, 1996).
\medskip

\item{} A.M.~Lorenzini: On the Betti numbers of points in projective 
space, Ph.D. thesis, Queen's University, Kingston, Ontario, 1987.
\medskip

\item{} A.M.~Lorenzini: The minimal resolution conjecture,
{\sl Journal of Algebra} {\bf 156} (1993), no. 1,  5--35.
\medskip

\item{} M.~Raynaud: Sections des fibr\'es vectoriels sur une courbe,
{\sl Bull. Soc. Math. France} {\bf 110} (1982), 103--125.
\medskip

\item{} F.-O.~Schreyer: Syzygies of canonical curves
with special pencils, Thesis, Brandeis University, 1983.
\medskip

\item{} C.~Walter: The minimal free resolution of the homogeneous ideal of
$s$ general points in $\P^4$, {\sl Math. Zeitschrift} {\bf 219} (1995),
no. 2, 231--234.
\medskip

}
\bigskip
\vbox{\noindent Author Addresses:
\smallskip
\noindent{David Eisenbud}\par
\noindent{Department of Mathematics, Brandeis University, 
Waltham MA 02254}\par
\noindent{eisenbud@math.brandeis.edu}
\medskip
\noindent{Sorin Popescu}\par
\noindent{Department of Mathematics, Brandeis University, 
Waltham MA 02254}\par
\smallskip
\noindent{current address:}\par
\noindent{Department of Mathematics, Columbia University, 
New  York, NY 10027}\par
\noindent{psorin@math.columbia.edu}\par
}

\end